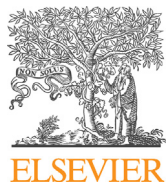
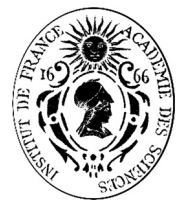

Algebra/Group theory

# On the $T^\circ$-slices of a finite group ☆

## Sur les $T^\circ$-tranches d'un groupe fini

Ibrahima Tounkara

*Laboratoire d'algèbre, de cryptologie, de géométrie algébrique et applications (LACGAA), Département de mathématiques et informatique, Faculté des sciences et techniques, Université Cheikh-Anta-Diop, BP 5005, Dakar, Senegal*

### A B S T R A C T

A slice $(G, S)$ of finite groups is a pair consisting of a finite group $G$ and a subgroup $S$ of $G$. In this paper, we show that some properties of finite groups extend to slices of finite groups. In particular, by analogy with $B$-groups, we introduce the notion of $T^\circ$-slice, and show that any slice of finite groups admits a largest quotient $T^\circ$-slice.

Crown Copyright © 2018 Published by Elsevier Masson SAS on behalf of Académie des sciences. All rights reserved.

### R É S U M É

Une tranche $(G, S)$ de groupes finis est un couple formé d'un groupe fini $G$ et d'un sous-groupe $S$ de $G$. Dans cet article, nous démontrons que certaines propriétés des groupes finis s'étendent aux tranches de groupes finis. En particulier, par analogie avec les $B$-groupes, nous introduisons la notion de $T^\circ$-tranche, et nous montrons que toute tranche de groupes finis admet un plus grand quotient qui soit une $T^\circ$-tranche.

## 1. Introduction

In this note, we extend to slices of finite groups some notions relative to finite groups. Throughout the paper, $G$ will denote a finite group.

In order to study the lattice of biset subfunctors of the Burnside functor $\mathbb{Q} \otimes B$, Serge Bouc (see [2]) studies the effect of the elementary biset operations on the primitive idempotents of the Burnside algebra $\mathbb{Q} \otimes B(G)$: these idempotents $e_H^G$ are indexed by the subgroups $H$ of $G$, up to conjugation and given (see [6], [8]) by

$$e_H^G = \frac{1}{|N_G(H)|} \sum_{K \leq H} |K| \mu(K, H) [G/K]$$

where $[G/K]$ denotes the isomorphism class of the $G$-set $G/K$, and $\mu$ is the Möbius function of the poset of subgroups of $G$.

Serge Bouc shows in particular that, for each normal subgroup $N$ of $G$, there is a rational number $m_{G,N}$ such that

$$\mathrm{Def}_{G/N}^{G} e_G^G = m_{G,N} e_{G/N}^{G/N}.$$

This constant $m_{G,N}$ is given by

$$m_{G,N} = \frac{1}{|G|} \sum_{\substack{X \leq G \\ XN=G}} |X| \mu(X,G).$$

This leads to the introduction of a special class of finite groups, called $B$-groups: a group $G$ is a $B$-group if $m_{G,N} = 0$ for any non-trivial normal subgroup $N$ of $G$.

When $N$ and $M$ are maximal subgroups of $G$ such that $m_{G,N} \neq 0$ and $m_{G,M} \neq 0$, then the quotients $G/N$ and $G/M$ are isomorphic. Such a quotient $G/N$ does not depend on $N$, up to isomorphism. It is denoted by $\beta(G)$. This quotient $\beta(G)$ of $G$ is a $B$-group, and moreover, any $B$-group $H$ that is a quotient of $G$ is actually a quotient of $\beta(G)$. Further results on this invariant $\beta(G)$ and $B$-groups can be found in [1] and [4].

Recall that a slice of finite groups is a pair $(T,S)$ consisting of a finite group $T$ and a subgroup $S$ of $T$. We say that the slice $(V,U)$ is a quotient of the slice $(T,S)$ if there exists a surjective group homomorphism $\varphi: T \to V$ such that $\varphi(S) = U$. When $\phi$ is an isomorphism, we say that $(T,S)$ and $(V,U)$ are isomorphic.

A slice of a finite group $G$ is a slice $(T,S)$ consisting of subgroups of $G$. The group $G$ acts by conjugation on the set $\Pi(G)$ of its slices.

In [3], the slice Burnside ring $\Xi(G)$ is introduced: it is a commutative ring, which has a $\mathbb{Z}$-basis $\langle T,S \rangle_G$ indexed by the conjugacy classes of slices $(T,S)$ of $G$. It is shown that $\mathbb{Q} \otimes \Xi(G)$ is a split semisimple $\mathbb{Q}$-algebra, whose primitive idempotents are also indexed by conjugacy classes of slices of $G$, and given by

$$\xi_{T,S}^G = \frac{1}{|N_G(T,S)|} \sum_{U \leq S \leq V \leq T} |U| \mu(U,S) \mu(V,T) \langle V,U \rangle_G.$$

It is also shown in [3] that the assignments $G \to \Xi(G)$ and $G \to \mathbb{Q} \otimes \Xi(G)$ are Green biset functors. On can try to describe the lattice of subfunctors of $\mathbb{Q} \otimes \Xi$, or its lattice of ideals.

In [7], it is shown that the idempotents $\xi_{G,S}^G$ play a similar role for $\mathbb{Q} \otimes \Xi$ as the idempotents $e_G^G$ do for $\mathbb{Q} \otimes B$: namely, for any $N \trianglelefteq G$, there exists a rational number $m_{G,S,N}$ such that

$$\mathrm{Def}_{G/N}^{G} \xi_{G,S}^G = m_{G,S,N} \xi_{G/N, SN/N}^{G/N}.$$

This constant $m_{G,S,N}$ is given by

$$m_{G,S,N} = \frac{|N_G(SN):SN|}{|N_G(S)|} \sum_{\substack{U \leq S \leq V \leq G \\ VN=G \\ UN=SN}} |U| \mu(U,S) \mu(V,G).$$

By Proposition 5.5 in [7] we have

$$m_{G,S,N} = \frac{|N_G(SN):SN|}{|N_G(S):S|} m_{S,S \cap N} m_{G,S,N}^{\circ}$$

where $m_{G,S,N}^{\circ} = \sum_{\substack{S \leq X \leq G \\ XN=G}} \mu(X,G)$. We observe that $m_{G,S,1} = m_{G,S,1}^{\circ} = 1$ for any slice $(G,S)$. Slices $(G,S)$ such that $m_{G,S,N} = 0$ for any non-trivial normal subgroup $N$ of $G$ have been called $T$-slices in [7]. By analogy, we say that a slice $(G,S)$ is a $T^{\circ}$-slice if $m_{G,S,N}^{\circ} = 0$ for any non-trivial normal subgroup $N$ of $G$.

As a first step towards a description of subfunctors of $\mathbb{Q} \otimes \Xi$, we state the following property of $T^{\circ}$-slices.

**Theorem 1.1.** *Let $(G,S)$ be a slice. Then there exists a slice $(H,U)$ such that*

- *$(H,U)$ is a $T^{\circ}$-slice and $(G,S) \twoheadrightarrow (H,U)$.*
- *if $(K,V)$ is a $T^{\circ}$-slice such that $(G,S) \twoheadrightarrow (K,V)$, then $(H,U) \twoheadrightarrow (K,V)$.*
  *Moreover, any two slices $(H,U)$ with these properties are isomorphic. We set $\tau^{\circ}(G,S) = (H,U)$.*

## 2. Proof of Theorem 1.1

We start with a proposition:

**Proposition 2.1.**

1. *Let $X$ be a subgroup of $G$, and $M$ be a normal subgroup of $G$. Then*

$$\mu(X, G) = \sum_{\substack{YM=G \\ Y \geq X \\ Y \cap M = X \cap M}} \mu(X, Y) \mu(Y, G).$$

2. *Let $S$ be a subgroup of $G$, and $M, N$ be normal subgroups of $G$. Then*

$$\mathrm{m}^\circ_{G,S,N} = \sum_{\substack{YN=YM=G \\ Y \geq S}} \mu(Y, G) \mathrm{m}^\circ_{G/M, SM/M, (Y \cap N)M/M}.$$

*In particular, if $N \geq M$, then $\mathrm{m}^\circ_{G,S,N} = \mathrm{m}^\circ_{G,S,M} \mathrm{m}^\circ_{G/M, SM/M, N/M}$.*

3. *If $M$ is a normal subgroup of $G$, maximal such that $\mathrm{m}^\circ_{G,S,M} \neq 0$, then $(G/M, SM/M)$ is a $T^\circ$-slice.*

**Proof.** 1. Let $[X, G]$ be the lattice of subgroups of $G$ containing $X$. In this lattice, a complement of $XM$ is a subgroup $Y$ of $G$ that contains $X$ such that $\langle Y, XM \rangle = G$ and $Y \cap XM = X$.
Since $X \leq Y$, the first condition is equivalent to $YM = G$, and the second one is equivalent to $X(Y \cap M) = X$, i.e. $Y \cap M = X \cap M$. Hence $|Y||M| = |G||X \cap M|$, so $|Y|$ depends only on $X$ and $M$. Hence there is no strict inclusion between complements of $XM$, and Crapo's formula (see [5]) gives the following result:

$$\mu(X, G) = \sum_{\substack{YM=G \\ Y \geq X \\ Y \cap M = X \cap M}} \mu(X, Y) \mu(Y, G)$$

for any normal subgroup $M$ of $G$.

2. By 1, we have

$$\mathrm{m}^\circ_{G,S,N} = \sum_{X,Y} \mu(X, Y) \mu(Y, G)$$

where the sum is over all pair $(X, Y)$ with

$$X \leq Y, \ XN = G, \ YM = G, \ Y \cap M = X \cap M, \ S \leq X \leq G$$

i.e.

$$YN = YM = G, \ X(Y \cap N) = Y, \ S \leq X \leq Y, \ X \geq (Y \cap M).$$

For a fixed $Y$, summing over $X$ is equivalent to summing over a subgroup $R \ (= X/(Y \cap M))$ of $Y/(Y \cap M)$ such that

$$R.(Y \cap N)(Y \cap M)/(Y \cap N) = Y/(Y \cap N), \ S.(Y \cap M)/(Y \cap M) \leq R \leq Y/(Y \cap M) \qquad (*)$$

Hence, the sum on $X$ is equal to

$$\sum_R \mu(R, Y/(Y \cap M)) = \mathrm{m}^\circ_{Y/(Y \cap M), S.(Y \cap M)/(Y \cap M), (Y \cap N)(Y \cap M)/(Y \cap N)} = \mathrm{m}^\circ_{G/M, SM/M, (Y \cap N)M/M},$$

where $R$ runs through subgroups of $Y/(Y \cap M)$ fulfilling Condition (*). The second equality here follows from the fact that since $YM = G$, there is a group isomorphism

$$Y/(Y \cap M) \cong G/M,$$

and this isomorphism sends $S(Y \cap M)/(Y \cap M)$ to $SM/M$ and $(Y \cap N)(Y \cap M)/(Y \cap M)$ to $(Y \cap N)M/M$.
Hence,

$$\mathrm{m}^\circ_{G,S,N} = \sum_{\substack{YN=YM=G \\ Y \geq S}} \mu(Y, G) \mathrm{m}^\circ_{G/M, SM/M, (Y \cap N)M/M}.$$

If $N \geq M$ and if $YM = G$, then $N = (Y \cap N)M$ and

$$m^{\circ}_{G,S,N} = \sum_{\substack{YM=G \\ Y \geq S}} \mu(Y,G) m^{\circ}_{G/M,SM/M,N/M} = m^{\circ}_{G,S,M} m^{\circ}_{G/M,SM/M,N/M}.$$

3. Straightforward. □

Let us give a proof of Theorem 1.1.

Let $S$ be a subgroup of $G$. Let moreover $N$ be a normal subgroup of $G$ such that $m^0_{G,S,N} \neq 0$, and let $(K,V)$ be a $T^{\circ}$-slice that is a quotient of $(G,S)$. Then there exists a surjective group homomorphism $\varphi : G \twoheadrightarrow K$ such that $\varphi(S) = V$. We denote by $M$ the kernel of $\varphi$, so $(K,V)$ is isomorphic to $(G/M, SM/M)$. Since $(G/M, SM/M)$ is a $T^{\circ}$-slice, Proposition 2.1 gives

$$m^{\circ}_{G,S,N} = \sum_{\substack{YN=YM=G \\ Y \cap N \leq M \\ Y \geq S}} \mu(Y,G).$$

Since $m^{\circ}_{G,S,N} \neq 0$, there exists a subgroup $Y$ of $G$ that satisfies:

$$YN = YM = G, \ (Y \cap N) \leq M, \ S \leq Y.$$

Then there exists a surjective homomorphism $Y/(Y \cap N) \twoheadrightarrow Y/(Y \cap M)$ that sends the group $S(Y \cap N)/(Y \cap N)$ to the group $S(Y \cap M)/(Y \cap M)$.

Since $Y/(Y \cap N) \cong G/N$, $Y/(Y \cap M) \cong G/M$, $S(Y \cap N)/(Y \cap N) \cong SN/N$, $S(Y \cap M)/(Y \cap M) \cong SM/M$, we have a surjective group homomorphism $G/N \twoheadrightarrow G/M$ sending $SN/N$ to $SM/M$. In other words,

$$(G/N, SN/N) \twoheadrightarrow (G/M, SM/M) \cong (K,V),$$

so $(K,V)$ is a quotient of $(G/N, SN/N)$.

If $M'$ is a normal subgroup of $G$, maximal such that $m^{\circ}_{G,S,M'} \neq 0$, then by Proposition 2.1 the slice $(H,U) := (G/M', SM'/M')$ is a $T^{\circ}$-slice, so

$$(G/N, SN/N) \twoheadrightarrow (H,U).$$

Moreover, since $m^{\circ}_{G,S,M'} \neq 0$, the slice $(K,V)$ is a quotient of $(G/M', SM'/M')$ i.e.

$$(H,U) \twoheadrightarrow (K,V).$$

Hence $(H,U)$ has the properties required in Theorem 1.1. If $(H',U')$ is another slice with these properties, then $(H,U)$ and $(H',U')$ are quotients of each other, so they are isomorphic. This completes the proof.

**Remark 2.2.** It is natural to ask if the analogue of Theorem 1.1 holds for $T$-slices instead of $T^{\circ}$-slices. Namely, for any slice $(G,S)$, does there exist a slice $\tau(G,S) = (H,U)$ such that

- $(H,U)$ is a $T$-slice and $(G,S) \twoheadrightarrow (H,U)$,
- if $(K,V)$ is a $T$-slice and $(G,S) \twoheadrightarrow (K,V)$, then $(H,U) \twoheadrightarrow (K,V)$?

The following example shows that the answer is no: consider the direct product $G = C_2 \times D_8$ of a group of order 2 and a dihedral group of order 8. Let $a$ denote the generator of $C_2$, and let $\{b,c\}$ be a set of generators of $D_8$, where $b$ has order 2 and $c$ has order 4. Set moreover $d = c^2$. Let $S$ be the subgroup of $G$ generated by $a$ and $b$. Thus $S \cong C_2 \times C_2$.

We denote by $N$ the subgroup of $G$ generated by $ad$, and $M$ the subgroup of $G$ generated by $d$.

Hence $|N| = |M| = 2$ and the subgroups $N$ and $M$ are central in $G$. We have $G/N \cong D_8$ and $G/M \cong (C_2)^3$.

The slices $(G/N, SN/N)$ and $(G/M, SM/M)$ are both quotients of $(G,S)$. On can check that they are both $T$-slices. If there exists a $T$-slice $(H,U)$ with the required properties, then in particular $H$ is a quotient of $G$, and both $G/N$ and $G/M$ are quotients of $H$. It follows that $H \cong G$, and then $(H,U)$ is isomorphic to $(G,S)$.

This is a contradiction, since $(G,S)$ is not a $T$-slice.

## Acknowledgements

The author is grateful to Serge Bouc (university of Picardie – Jules-Verne) and Prof. Oumar Diankha (university of Cheikh Anta Diop) for their comments and guidance.

364I. Tounkara / C. R. Acad. Sci. Paris, Ser. I 356 (2018) 360–364